\documentclass[10pt]{amsart}

\usepackage{amssymb}
\usepackage{amsthm}

\theoremstyle{plain}
\newtheorem{theorem}{Theorem}[section]
\newtheorem{lemma}[theorem]{Lemma}

\newtheorem{claim}[theorem]{Claim}

\theoremstyle{definition}
\newtheorem{definition}{Definition}[section]

\newtheorem{question}{Question}[section]

\theoremstyle{remark}
\newtheorem*{remark}{Remark}

\newtheorem{case}{Case}

\newcommand{\restr}{\mathrel{\mbox{\raisebox{.5mm}{$\upharpoonright $}}}}
\newcommand{\D}{{\mathcal D}}
\newcommand{\F}{{\mathcal F}}
\renewcommand{\phi}{{\varphi}}

\newif\ifpdf
\ifx\pdfoutput\undefined
\pdffalse 
\else
\pdfoutput=1 
\pdftrue
\fi

\ifpdf
\fi

\begin{document}

\title{Turing Incomparability in Scott Sets}

\author{Anton\'in Ku\v{c}era}
\address{Department of Theoretical Computer Science and Mathematical Logic\\
  Faculty of Mathematics and Physics\\
  Charles University\\
  Malostransk\'{e} n\'{a}m. 25, 118 00 Praha 1\\
  Czech Republic}
\email{kucera@ksi.mff.cuni.cz}

\thanks{Ku\v{c}era was partially supported by the Research project
  of the Ministry of Education of the Czech Republic MSM0021620838}

\author{Theodore A. Slaman}
\address{Department of Mathematics\\
  The University of California, Berkeley\\
  Berkeley, CA 94720-3840 USA}
\email{slaman@math.berkeley.edu}
\thanks{Slaman was partially supported by NSF grant DMS-0501167.}


\begin{abstract}
  For every Scott set $\mathcal F$ and every non-recursive set $X$ in
  $\mathcal F$, there is a $Y \in \mathcal F$ such that $X$ and $Y$
  are Turing incomparable.
\end{abstract}

\maketitle

\section{Introduction}

H. Friedman and A. McAllister posed the question whether for every
non-recursive set $X$ of a Scott set $\mathcal F$ there is a $Y \in
\mathcal F$ such that $X$ and $Y$ are Turing incomparable (see Problem
3.2 and also Problem 3.3 in Cenzer and Jockusch
\cite{Cenzer.Jockusch:2000}).  We present a positive solution to the
question, using recent results in the area of algorithmic randomness
and also results on $\Pi^0_1$ classes.

\subsection{Background and Notation}  We begin by discussing the
background of the Friedman-McAllister question.  We then review some
basic definitions and establish our notational conventions.

\subsubsection{Background}

We let $2^\omega$ denote the set of infinite binary sequences.  One
can equivalently think of $2^\omega$ as the Cantor set.  A finite
binary sequence $\sigma$ determines an open neighborhood in $2^\omega$
by taking the set of all infinite extensions of $\sigma$.  A 
binary tree $T$ determines a closed subset of $2^\omega$ by taking the
complement of the union of open neighborhoods given by the elements of
$T$ which have no extensions in $T$.  

As is well-known, the Cantor set is a canonical example of a compact
set.  This fact translates to binary trees in the form of
K\"onig's Lemma that every infinite binary tree has an infinite path.
However, the proof of K\"onig's Lemma is not computational.  Not every
infinite recursive binary tree has an infinite recursive path.  Thus, the set of
recursive reals does not verify the compactness of $2^\omega$ with
respect to recursive closed sets (also called $\Pi^0_1$ classes) .  In
order to study the consequences of compactness, we need richer subsets
of $2^\omega$.

\begin{definition}
  A Scott set is a nonempty set $\mathcal F \subseteq 2^\omega$ such
  that if $T \subseteq 2^{<\omega}$ is an infinite tree
  recursive in a finite join of elements of $\mathcal F$, then
  there is an infinite path through $T$ in $\mathcal F$.
\end{definition}

Scott \cite{Scott:1962} proved that the sets representable in a
complete extension of Peano Arithmetic form a Scott set.  Scott sets
occur naturally in the study of models of arithmetic.  They are the
$\omega$-models of $WKL_0$, an axiomatic treatment of compactness.  In
other words, they are the models of $WKL_0$ in which the natural
numbers are standard.

One can test the power of compactness arguments by examining what is
true in every Scott set.  A natural family of questions comes from
considering the Turing degrees represented in an arbitrary Scott set.
For example, it is not merely the case that every Scott set has
non-recursive elements.  Jockusch and Soare~\cite{Jockusch.Soare:1972}
showed that for any Scott set $S$ and any finite partial order $P$,
there are elements of $S$ whose Turing degrees are ordered
isomorphically to $P$.  Thus, the existential theory of the Turing
degrees of a Scott set is rich and completely determined.

The existential-universal theory of the degrees represented in an
arbitrary Scott set is more complex and not at all understood.
Groszek and Slaman \cite{Groszek.Slaman:1997} states that every Scott
set has an element of minimal Turing degree, namely a degree $m$ such
that $m$ has no non-trivial element strictly below it.
Friedman-McAllister question is universal-existential about an
arbitrary Scott set, for every degree $x$ is there a degree $y$ which
is Turing incomparable with $x$.  In other words, given a Turing degree $x$,
can one use a compactness argument to construct a $y$ which is Turing
incomparable with $x$?

The typical immediate reaction to the question is that there should
always be such a $y$ and that it should be routine to exhibit an
$x$-recursive tree such that every infinite path in that tree has
Turing degree incomparable with $x$.  This is not the case.  For
example, Ku\v{c}era \cite{Kucera:1988} showed that there is a Scott
set $S$ and a non-recursive degree $x$ from $S$ such that $x$ is
recursive in all complete extensions of Peano Arithmetic that appear
in $S$.  Building $y$'s incomparable with that $x$ cannot be
accomplished using complete extensions of Peano Arithmetic.  Similar
obstacles appear when one attempts to find $y$'s incomparable to $x$
by means of other familiar $\Pi^0_1$ classes.

Even so, for every Scott set $S$ and every non-recursive $x$
represented in $S$, there is a $y$ in $S$ which is Turing incomparable
with $x$.  Given $x$, our construction of $y$ is not uniform.  If
possible, we find $y$ by taking a sequence which is 1-random relative
to $x$.  If that fails (the non-uniformity), then we apply recent
results in the theory of algorithmic randomness to build a recursive
tree whose infinite paths are Turing incomparable with $x$.

\subsubsection{Definitions and Notation}

Our computability-theoretic notation generally follows Soare
\cite{Soare:1987} and Odifreddi \cite{Odifreddi:1989,Odifreddi:1999}.
An introduction to algorithmic randomness can be found in Li and
Vit\'anyi \cite{Li.Vitanyi:1997}.  A short survey of it is also given
in Ambos-Spies and Ku\v{c}era \cite{Ambos-Spies.Kucera:2000}.  Much
deeper insight into the subject of algorithmic randomness can be found
in a forthcoming book of Downey and Hirschfeldt \cite{Downey.Hirschfeldt:nd}, a good
survey is also in Downey, Hirschfeldt, Nies and Terwijn \cite{Downey.Hirschfeldt.ea:nd}.

We refer to the elements of $2^{\omega}$ as sets or infinite binary
sequences.  We denote the collection of strings, i.e. finite initial
segments of sets, by $2^{<\omega}$.  The length of a string $\sigma$
is denoted by $|\sigma|$, for a set $X$, we denote the string
consisting of the first $n$ bits of $X$ by $X \restr n$ and we use
similar notation $\sigma \restr n$ for strings $\sigma$ of length
$\geq n$.  We let $\sigma * \tau$ denote the concatenation of $\sigma$
and $\tau$ and let $\sigma * Y$ denote the concatenation of $\sigma$
and (infinite binary sequence) $Y$.  We write $\sigma \prec X$ to
indicate $X\restr|\sigma|=\sigma$.  If $\sigma \in 2^{<\omega}$, then
$[\sigma]$ denotes $\{X \in 2^\omega : \sigma \prec X \}$.

A $\Sigma^0_1$ class is a collection of sets that can be effectively
enumerated.  Such a class can be represented as $\bigcup_{\sigma \in
  W}[\sigma]$ for some (prefix-free) recursively enumerable (r.e.) set
of strings $W$.  The complements of $\Sigma^0_1$ classes are called
$\Pi^0_1$ classes. Any $\Pi^0_1$ class can be represented by the class of
all infinite paths through some recursive tree. We use also
relativized versions, i.e.  $\Sigma^{0, X}_1$ classes and $\Pi^{0,
  X}_1$ classes.  $\Pi^0_1$ classes play an important role in logic,
in subsystems of second-order arithmetic, and also in algorithmic
randomness.  By the relativized tree representation of $\Pi^0_1$ classes, if
$\mathcal F$ is a Scott set, $X\in\mathcal F$, and $P$ is a nonempty
$\Pi^{0,X}_1$ class, then $\mathcal F$ includes an element of $P$.

\begin{definition}
  Let $X$ be a set. A Martin-L\"{o}f test relative to $X$ is a
  uniformly r.e. in $X$ sequence of $\Sigma^{0, X}_1$ classes
  $\{U^X_n\}$ such that $\mu(U^X_n) \leq 2^{-n}$, where $\mu$ denotes
  the standard measure on $2^\omega$.  Then any subclass of
  $\bigcap_{n\in\omega}U^X_n$ is called a Martin-L\"{o}f null set
  relative to $X$.  If $X = \emptyset$, we simply speak of
  Martin-L\"{o}f test and Martin-L\"{o}f null set.  A set $R$ is
  Martin-L\"{o}f random relative to $X$, or $1$-random relative to
  $X$, if $\{R\}$ is not Martin-L\"{o}f null relative to $X$. If $X =
  \emptyset,$ we speak of $1$-randomness.
\end{definition}

Martin-L\"{o}f proved that there is a universal Martin-L\"{o}f test,
$\{U_n\}$, such that for all $R$, $R$ is 1-random if and only if
$R\not\in\bigcap_{n\in\omega}U_n.$ Similarly, there is a universal
Martin-L\"{o}f test relative to $X,$ $\{U^X_n\}$ (uniformly in $X$).

We will use $K$ to denote prefix-free Kolmogorov complexity.  See Li
and Vit\'anyi \cite{Li.Vitanyi:1997} for details.  The version of
Kolmogorov complexity relativized to a set $X$ is denoted by $K^X$.
Schnorr \cite{Schnorr:1971*b} proved that a set is $1$-random if and only if for
all $n$, $K(A \restr n) \geq n + O(1)$ .  There are several notions of
computational weakness related to $1$-randomness.  They are summarized
in the following definition.

\begin{definition}
  \begin{enumerate}
  \item $\mathcal L$ denotes the class of sets which are low for
    $1$-randomness, i.e sets $A$ such that every $1$-random set is
    also $1$-random relative to $A$. 
  \item ${\mathcal K}$ denotes the class of $K$-trivial sets, i.e. the
    class of sets $A$ such that for all $n,$ $K(A \restr n) \leq
    K(0^n) + O(1)$, where $0^n$ denotes the string of $n$ zeros.
  \item $\mathcal M$ denotes the class of sets that are low for $K$,
    i.e sets $A$ such that for all $\sigma$, $K(\sigma) \leq
    K^A(\sigma) + O(1).$
  \item  A set $A$ is a basis for $1$-randomness if $A \leq_T Z$ for
    some $Z$ such that $Z$ is $1$-random relative to $A$. The
    collection of such sets is denoted by $\mathcal B$. 
  \end{enumerate}
\end{definition}

Nies \cite{Nies:2005} proved that $\mathcal L$ = $\mathcal M$,
Hirschfeldt and Nies, see \cite{Nies:2005}, proved that $\mathcal K$ =
$\mathcal M$, and Hirschfeldt, Nies and Stephan
\cite{Hirschfeldt.Nies.ea:nd} proved that $\mathcal B = \mathcal K$.
Thus, all these four classes are equal and we have, remarkably, four
different characterizations of the same class.  That is, $\mathcal L =
\mathcal K = \mathcal M = \mathcal B$.

Chaitin \cite{Chaitin:1977} proved that $K$-trivials are $\Delta^0_2$.
Further, by a result of Ku\v{c}era \cite{Kucera:1993} low for
$1$-random sets are $GL_1$ and, thus, all these sets are, in fact,
low.  The lowness of these sets also follows from some recent results
on this class of sets, see \cite{Nies:2005}.

\section{The main result} 

In this section, we present a solution to the Friedman--McAllister
question.

\begin{theorem}\label{2.1}
  For any Scott set $\mathcal S$ and any non-recursive set
  $X \in \mathcal S$, there is a $Y \in \mathcal S$ such that
  $X$ and $Y$ are Turing incomparable.
\end{theorem}

Theorem~\ref{2.1} is a consequence of the stronger Claim~\ref{2.2}.

\begin{claim}\label{2.2}
  For every non-recursive set $X$ there is a nonempty $\Pi^{0, X}_1$ class
  $\mathcal P \subseteq 2^{\omega}$ such that  every element of
  $\mathcal P$  is Turing incomparable with $X$.
\end{claim}

As we will see, we can do even better for $K$-trivial sets $X$.
Namely, we can replace the $\Pi^{0, X}_1$ class mentioned in the claim
by a (non-relativized) $\Pi^0_1$ class.

\begin{proof}[Proof of Claim~\ref{2.2}]
  
  We split the proof of the claim into two cases.

  \begin{case}
    $X$ is not a basis for $1$-randomness.
  \end{case}

  Let $T_1$ be a tree recursive in $X$ such that any infinite path in
  $T_1$ is $1$-random relative to $X$.  We can take e.g. a tree
  recursive in $X$ for which the collection of all infinite paths is the
  $\Pi^{0,X}_1$ class which is the the complement of $U^X_0$, the first
  class appearing in a universal Martin-L\"of test relative to $X$.
  Clearly, any infinite path through $T_1$ is Turing incomparable with
  $X$.

  \begin{case}
    $X$ is a basis for $1$-randomness.
  \end{case}

  As we described above, such an $X$ is $K$-trivial, low for
  $1$-randomness, low for $K$, and $\Delta^0_2$.  We use these
  properties to construct a recursive tree $T_2$ such that any infinite
  path through $T_2$ is Turing incomparable with $X$. 

  Also we not only avoid a lower cone of sets recursive in $X$, but even
  avoid all the class of sets which are bases for $1$-randomness. Since
  the class of these sets is closed downwards, we obviously get a
  stronger property.

  Since the class $\mathcal B$ is equal to $\mathcal L$ and also to
  $\mathcal M$, we can equivalently work with any of these
  characterizations. 

  Thus, to handle Case~2, it is sufficient to prove the following lemma.

  \begin{lemma}\label{2.3}
    Let $X$ be a non-recursive $\Delta^0_2$ set.  Then there is a
    recursive tree $T_2$ such that for any infinite path $Y$ in $T_2$ we
    have $Y \not \geq_T X$ and $Y$ is not low for $1$-randomness (and,
    thus, $Y \not \leq_T X$ if $X \in \mathcal L$).
  \end{lemma}

  \begin{proof}[Proof of Lemma~\ref{2.3}]
    
The proof is by a finite injury priority argument. We build the tree
$T_2$ by stages.  At stage $s+1$, we terminate a string by not
extending it to any string of length $s+1$ in $T_2$.
    
We will describe the strategies and leave the rest to the reader.

The strategies have the following general pattern.  Each strategy
starts to work at a given string $\sigma \in T_2$, it acts only
finitely often, and it yields as its outcome a nonempty finite
collection $Q$ consisting of strings of the same length.  Some
strategies (called $\mathcal L$-strategies) and their outcomes depend
not only on $\sigma$ itself, but also on how $\sigma$ arises as a
concatenation of strings belonging to outcomes of previous strategies.
Further, for each $\alpha \in Q$, the string $\sigma*\alpha$ together
with a recursive tree of all strings extending $\sigma*\alpha$ are
left for the next strategy.  By producing its outcome, each strategy
satisfies some particular requirement as explained later.

\textit{Avoiding an upper cone above $X$.\ \ } Let $\sigma \in T_2$ be
given. Suppose that $\Phi$ is a recursive Turing functional. We act to
ensure that for every infinite path $A$ in $T_2$, if $A$ extends
$\sigma$ then $\Phi(A)\neq X$.

We use the fact that $X$ is $\Delta^0_2$ to adapt the Sacks
Preservation Strategy \cite{Sacks:1963*b}.  We monitor the maximum length of
agreement between $\Phi(\tau)$, for $\tau$ extending $\sigma$, and the
current approximation to $X$.  If at stage $s+1$, we see a string
$\eta$ of length $s$ on $T_2$ for which this maximum has gone higher
than ever before, we take the least such $\eta$ and we terminate all
extensions of $\sigma$ except for $\eta$.

If this were to occur infinitely often above $\sigma$, then $X$ would
be recursive. Compute $X(n)$ by finding the first stage where the
maximum length of agreement between $\Phi(\tau)$, for some $\tau$
extending $\sigma$, and the approximation to $X$ was greater than $n$.
Since the length of agreement increases infinitely often, the
approximation to $X$ returns to this value infinitely often.  But
then, since the approximation converges to the value of $X$, the value
at the stage we found must be the true value.  This is a contradiction
to $X$'s being non-recursive.

So, this strategy acts finitely often and satisfies the requirement.

Observe, that the strategy yields as its output just one string
$\alpha$, where $\alpha$ is the string at which we last terminate all
extensions of $\sigma$ which are incompatible with $\alpha$.
Otherwise, $\alpha$ is the empty string, if we never do so.

\textit{Avoiding the class of sets low for $1$-randomness.\ \ } We
will refer to our strategies to avoid the class of sets low for
$1$-randomness as $\mathcal L$-strategies.  We will begin by
explaining the general idea behind these strategies.
    
Recall, a set $X$ is low for $1$-randomness if and only if every
1-random set $Z$ is also 1-random relative to $X$.  In the case that
$X$ is low for $1$-randomness, we can ensure that $X$ does not compute
any path in $T_2$ by ensuring that each path $Y$ in $T_2$ is not low for
1-randomness, i.e. $Y\not \in\mathcal L$.

We will ensure that a path $Y$ in $T_2$ is not in $\mathcal L$, by
embedding large intervals of some $1$-random set $Z$ into it. In this way, we
can recover the 1-random set $Z$ recursively from $Y$ and $\emptyset'$
and ensure that $Y$ can enumerate a Martin-L\"{o}f test (relative to
$Y$) which shows that $Z$ is not $1$-random in $Y$.
    
An infinite path $Y$ in $T_2$ can be viewed as an
infinite concatenation of strings $\alpha_0 * \alpha_1 * \alpha_2 *
\dots $ , where each $\alpha_i$ is that uniquely determined string
compatible with $Y$, which belongs to the outcome of a strategy that
started at $\alpha_0 * \dots * \alpha_{i-1}$ (where $\alpha_{-1}$ is
the empty string).  Let us note that due to a standard finite injury
priority argument, such a sequence $\{\alpha_i \}_{i \in \omega}$ can
be found recursively in $\emptyset'$ and $Y$.  Then let $Z_Y$ denote
the set obtained as an infinite concatenation of strings $\alpha_{i_0}
* \alpha_{i_1} * \alpha_{i_2} * \ldots $, where $\{i_j \}_{j \in
  \omega}$ is a recursive increasing sequence of indices of those
strings $\alpha_i$'s which belong to outcomes of $\mathcal
L$-strategies, i.e. of those $i$'s for which a strategy that started
at $\alpha_0 * \alpha_1 * \ldots * \alpha_{i-1}$ was an $\mathcal
L$-strategy.  The general goal of $\mathcal L$-strategies is to ensure
that for any infinite path $Y$ in $T_2$, $Z_Y$ is $1$-random but is
not $1$-random relative to $Y$.

To guarantee that $Z_Y$ is not $1$-random in $Y$, we have to satisfy
for all $e$ the requirement $Z_Y \in U^Y_e$, where \{$U^Y_e\}_{e
  \in\omega}$ is a universal Martin-L\"{o}f test relative to $Y$
(uniformly in $Y$).  As is standard, we may let $U^Y_e$ be $\bigcup_k
V^{k,Y}_{e+k+1}$, where $\{V^{k,Y}_i \}_{i,k \in \omega}$ is a
uniformly r.e. in $Y$ sequence of all Martin-L\"{o}f
tests relative to $Y$ (uniformly in $Y$), and $\{V^{k,Y}_i \}_{i
  \in\omega}$ is the $k$-th test.

To guarantee that $Z_Y$ is $1$-random, we fix a $\Pi^0_1$ class
$\overline{U}_0$ of 1-random sequences and fix a recursive tree $T^*$
such that the infinite paths in $T^*$ are exactly the members of
$\overline{U}_0$.  We will ensure that each initial segment of $Z_Y$
extends to an element of $\overline{U}_0$.
    
Suppose now, that $\sigma \in T_2$ and $e$ are given, and let
$\alpha_0, \ldots , \alpha_k$ be the strings for which $\sigma =
\alpha_0 * \alpha_1 * \dots * \alpha _k $, where each $\alpha_i$
belongs to an outcome of the strategy that started at $\alpha_0 *
\alpha_1 * \dots * \alpha_{i-1}$.  Let $\tau_\sigma$ be the string
$\alpha_{i_0} * \alpha_{i_1} * \ldots * \alpha_{i_j}$ , where $i_0,
i_1, \ldots , i_j$ are indices (in increasing order) of those
$\alpha_i$'s which belong to outcomes of $\mathcal L$-strategies, i.e.
of those $i$'s such that the strategy that started at $\alpha_0 *
\alpha_1 * \ldots * \alpha_{i-1}$ was an $\mathcal L$-strategy.
Roughly speaking, $\tau_\sigma$ is the finite sequence already
embedded into $\sigma$ which can be extended to an infinite path in
$T^*$.

Observe, that for any set $Z$ the set $\tau_\sigma*Z$ is not
$1$-random relative to the set $\sigma*Z$.  In fact, $\tau_\sigma*Z$
is recursive in $\sigma*Z$ as it is obtained by appending all but
finitely much of $\sigma*Z$ to $\tau_\sigma$.  Thus, $\tau_\sigma*Z\in
U^{\sigma*Z}_e$.  If we had no other requirements to satisfy, we could
restrict the infinite extensions of $\sigma$ in $T_2$ to those of the
form $\sigma*Z$ for which $\tau_\sigma*Z$ is $1$-random (using the
$\Pi^0_1$ class $\overline{U}_0$).  But, of course, we have to satisfy
our requirement in a finitary way to leave space for the cone avoiding
strategies.  The idea here is to make any infinite path through $T_2$
extending $\sigma$ locally $1$-random.

Thus, we design our tree so that enough of such a $Z$ is
embedded in the extensions of $\sigma$ to ensure that $[\tau_\sigma*(Z
\restr i)] \subseteq U^{\sigma*Z}_e $ for some $i$.  The crucial thing
here is that we can accomplish this objective in a finite way.  That
is, we can effectively compute (from $\sigma$, $\tau_\sigma$ and $e$)
an $i$ such that $[\tau_\sigma*(Z \restr i)] \subseteq U^{\sigma*Z}_e$
for all sets $Z$.  Intuitively, for any $Z$, $\tau_\sigma*Z$ is not
1-random relative to $\sigma*Z$ and we can calculate how long it takes for
$\sigma*Z$ to recognize the failure of relative 1-randomness.

We give this calculation in detail.  Given $\sigma$ and $\tau_\sigma$
find a Martin-L\"{o}f test relative to $X$ (uniformly in $X$) $\{B^X_j
\}_{j \in \omega}$ with index $b$ such that $B^X_j =
[(\tau_\sigma*X^*) \restr j]$, where $X^*$ is the set for which $X =
(X \restr |\sigma|) *X^*$.  Then we obviously have $B^{\sigma*Z}_j =
[(\tau_\sigma*Z) \restr j]$, for any set $Z$.  By the construction of
$\{U^X_e\}_{e \in\omega}$, the universal Martin-L\"{o}f test (relative
to $X$), and since $b$ is an index of the test $\{B^X_j \}_{j \in
  \omega}$, we have $B^{\sigma*Z}_{e+b+1} \subseteq U^{\sigma*Z}_e$
for all sets $Z$.  It follows that $[\tau_\sigma*(Z \restr i)]
\subseteq U^{\sigma*Z}_e$ for all sets $Z$ and $i$ such that
$|\tau_\sigma| + i \geq e+b+1$.  Our calculation chooses the least such
$i$.

It only remains to put a restriction on $T_2$ to ensure that $\sigma *
\alpha$ is extendable to an infinite path in $T_2$ for strings
$\alpha$ of length $i$ if and only if $\tau_\sigma * \alpha$ is
extendable to an infinite path in $T^*$, the recursive tree whose
infinite paths are exactly the elements of $\overline{U}_0$ and hence
are 1-random.

The strategy, given $\sigma \in T_2$ and $e$, where $\sigma = \alpha_0
* \alpha_1 * \ldots * \alpha_k$ with properties of $\alpha_i$'s
described above, is precisely as follows.  Find the corresponding
string $\tau_\sigma$, then compute an $i$ such that $[\tau_\sigma*(Z
\restr i)] \subseteq U^{\sigma*Z}_e$ for all sets $Z$.  Now for each
$\beta$ such that $|\beta| \geq i$, we terminate the string
$\sigma*\beta$ in $T_2$ if $\tau_\sigma * (\beta \restr i)$ is not
extendable to a string of length $|\sigma * \beta|$ in the recursive
tree $T^*$ which represents the $\Pi^0_1$ class $\overline{U}_0$.  This
strategy acts only finitely often and eventually reaches its goal.
Observe that the strategy yields a finite collection of strings of the
same length $Q$, where all requests on $Q$ are satisfied.
  \end{proof}

  With Lemma~\ref{2.3}, we have completed the proof of
  Claim~\ref{2.2}.\end{proof}

\subsection{An $\mathcal M$ variation}

Since $\mathcal L$ and $\mathcal M$ are equal, we
can equivalently use the characterization of $\mathcal M$ to design
our strategy to handle Case~2, namely, our strategy for 
avoiding the class of sets which are bases for $1$-randomness. 
For the convenience of the reader we present
also a variant of a strategy expressed in terms of $\mathcal M$.

Given a $\sigma$, we want to ensure that each infinite path $Y$
extending $\sigma$ in $T_2$ can give a shorter description of some
string $\tau$ than any description possible without $Y$.  Let $c$ be
the amount that we want to shorten the description.  We will compute
an $m$ (see below) and we want to ensure
\[
K(Y \restr m) - c \geq K^Y(Y\restr m).
\]

We choose $m$ much larger than $K(\sigma)$ and $c$. The maximum of 10
and $2^{|\sigma|+K(\sigma)+c+d}$ is big enough (where a constant $d$
is explained below).  For each string $\tau$ extending $\sigma$ of
length $m$ let $\tau^*$ denote the string of length $m- |\sigma|$ for
which $\tau = \sigma * \tau^*$.  It is easy to see that $K^Y(Y \restr
m)$ is less than or equal to $2 \log(m)$, since $Y$ can describe its
first $m$ values using the description of $m$.  (As a caveat, this
bound may only apply to sufficiently large $m$ because of the fixed
cost of interpreting binary representations. This is fixed data and we
can assume that $m$ is large enough for the upper bound to apply).

On the other hand, $K(\tau^*) \leq K(\tau) + K(\sigma) + d$ for some
constant $d$ independent of a choice of $\tau$.  By terminating
strings with shorter descriptions, we can ensure that for each $Y$
extending $\sigma$ in $T_2$, the $\tau^*$ of $\tau = Y \restr m$
satisfies $K(\tau^*) \geq m - |\sigma|$. This is similar to making a
recursive tree of $1$-random sets, but here we are making any path in 
$T_2$ (merely) locally $1$-random to ensure that infinite paths in
$T_2$ are not low for $K$.

We can now calculate:
\[
K(\tau) \geq K(\tau^*) - K(\sigma) - d,
\]
and substituting for $K(\tau^*)$, 
\[
K(\tau) \geq (m - |\sigma|) - K(\sigma) - d.
\]
Since $K^Y(Y \restr m) \leq 2\log(m)$, it is sufficient 
to ensure that 
\[
m - |\sigma| -K(\sigma) -d - c \geq 2\log(m),
\]
or, equivalently, 
\[
m \geq 2\log(m) + |\sigma| + K(\sigma) + c + d.
\]
If $m \geq 2^{|\sigma| + K(\sigma)+c+d}$, then it is sufficient to ensure 
$m \geq 3\log(m)$. This holds if $m$ is greater than 10.

So, the strategy working above $\sigma$ to ensure that $Y$ is not low
for $K$ reserves the collection of extensions of $\sigma$ of length
$m$ and at most half of them are eventually stopped to be extendable
to an infinite path in $T_2$. So, it satisfies its requirement and
acts only finitely often.

\begin{remark}
  We have not addressed the question whether it is provable in the
  subsystem of second order arithmetic $WKL_0$ that for every
  non-recursive set $X$ there is a non-recursive set $Y$ which is Turing
  incomparable with $X$ (see Problem 3.2. part 1 in Cenzer and
  Jockusch \cite{Cenzer.Jockusch:2000}).
\end{remark}

\section{An open problem}

Suppose that $\F$ is a Scott set and let $\D^\F$ denote the partial
order of the Turing degrees which are represented by elements of $\F$.
According to Theorem~\ref{2.1},
\[
\D^\F\models \forall d>0\exists x(d\not\geq_T x \text{ and }x\not\geq_T d).
\]
The dual theorem of Groszek and Slaman \cite{Groszek.Slaman:1997} states
that every Scott set has an element of minimal Turing degree:
\begin{equation}
  \label{eq:minimal}
  \D^\F\models \exists d>0\forall x \neg(d>_T x \text{ and }x>_T 0).
\end{equation}

Together, these results are sufficient to determine for any sentence in
the language of partial orders of the form $\forall d\exists
x\phi(d,x)$, where $\phi$ is quantifier-free, whether that sentence
holds in $\D^\F$.  Further, such sentences hold in $\D^\F$ if and only
if they hold in $\D$, the Turing degrees of all sets.

By Lerman~\cite{Lerman:1983} and Shore~\cite{Shore:1978}, the general
$\forall\exists$-theory of $\D$ is decidable.  The proof of
decidability rests on two technical results.  The first is a general
extension theorem due to Sacks~\cite{Sacks:1961*b} which, like
Theorem~\ref{2.1}, constructs degrees $x$ incomparable to or above
given ones $d$.  The second is Lerman's~\cite{Lerman:1971} theorem
that every finite lattice is isomorphic to an initial segment of the
Turing degrees which, like Formula~(\ref{eq:minimal}), produces
degrees $d$ which limit the possible types of $x$'s which are below
$d$.

Superficially, Theorem~\ref{2.1} and the existence of minimal degrees
suggest that the $\forall\exists$-theory of $\D^\F$ resembles that of
$\D$.  However, the actual proofs are quite different, and we are left
with the following question.

\begin{question}
  Suppose that $\F$ is a Scott set.  Is $\D^\F$
  $\forall\exists$-elementarily equivalent to $\D$?
\end{question}


\begin{thebibliography}{10}

\bibitem{Ambos-Spies.Kucera:2000}
Klaus Ambos-Spies and Anton{\'{\i}}n Ku{\v{c}}era.
\newblock Randomness in computability theory.
\newblock In {\em Computability theory and its applications (Boulder, CO,
  1999)}, volume 257 of {\em Contemp. Math.}, pages 1--14. Amer. Math. Soc.,
  Providence, RI, 2000.

\bibitem{Cenzer.Jockusch:2000}
Douglas Cenzer and Carl~G. Jockusch, Jr.
\newblock {$\Pi\sb 1\sp 0$} classes---structure and applications.
\newblock In {\em Computability theory and its applications (Boulder, CO,
  1999)}, volume 257 of {\em Contemp. Math.}, pages 39--59. Amer. Math. Soc.,
  Providence, RI, 2000.

\bibitem{Chaitin:1977}
G.~J. Chaitin.
\newblock Algorithmic information theory.
\newblock {\em IBM J. Res. Develop}, 21:350--359, 496, 1977.

\bibitem{Downey.Hirschfeldt:nd}
R.~Downey and D.~R. Hirschfeldt.
\newblock Algorithmic randomness and complexity.
\newblock to appear.

\bibitem{Downey.Hirschfeldt.ea:nd}
R.~Downey, D.~R. Hirschfeldt, A.~Nies, and S.~A. Terwijn.
\newblock Calibrating randomness.
\newblock to appear.

\bibitem{Groszek.Slaman:1997}
Marcia~J. Groszek and Theodore~A. Slaman.
\newblock ${\Pi}\sp 0\sb 1$ classes and minimal degrees.
\newblock {\em Ann. Pure Appl. Logic}, 87(2):117--144, 1997.
\newblock Logic Colloquium '95 Haifa.

\bibitem{Hirschfeldt.Nies.ea:nd}
D.~R. Hirschfeldt, A.~Nies, and F.~Stephan.
\newblock Using random sets as oracles.
\newblock to appear.

\bibitem{Jockusch.Soare:1972}
Carl~G. Jockusch, Jr. and Robert~I. Soare.
\newblock {$\Pi_1^0$} classes and degrees of theories.
\newblock {\em Trans. Amer. Math. Soc.}, 173:33--56, 1972.

\bibitem{Kucera:1988}
Anton{\'{\i}}n Ku{\v{c}}era.
\newblock On the role of {${\bf 0}'$} in recursion theory.
\newblock In {\em Logic colloquium '86 (Hull, 1986)}, volume 124 of {\em Stud.
  Logic Found. Math.}, pages 133--141. North-Holland, Amsterdam, 1988.

\bibitem{Kucera:1993}
Anton{\'{\i}}n Ku{\v{c}}era.
\newblock On relative randomness.
\newblock {\em Ann. Pure Appl. Logic}, 63(1):61--67, 1993.
\newblock 9th International Congress of Logic, Methodology and Philosophy of
  Science (Uppsala, 1991).

\bibitem{Lerman:1971}
M.~Lerman.
\newblock Initial segments of the degrees of unsolvability.
\newblock {\em Ann. of Math.}, 93:311--389, 1971.

\bibitem{Lerman:1983}
M.~Lerman.
\newblock {\em Degrees of Unsolvability}.
\newblock Perspectives in Mathematical Logic. Springer--Verlag, Heidelberg,
  1983.
\newblock 307 pages.

\bibitem{Li.Vitanyi:1997}
Ming Li and Paul Vit{\'a}nyi.
\newblock {\em An introduction to {K}olmogorov complexity and its
  applications}.
\newblock Graduate Texts in Computer Science. Springer-Verlag, New York, second
  edition, 1997.

\bibitem{Nies:2005}
Andr{\'e} Nies.
\newblock Lowness properties and randomness.
\newblock {\em Adv. Math.}, 197(1):274--305, 2005.

\bibitem{Odifreddi:1989}
Piergiorgio Odifreddi.
\newblock {\em Classical recursion theory}, volume 125 of {\em Studies in Logic
  and the Foundations of Mathematics}.
\newblock North-Holland Publishing Co., Amsterdam, 1989.
\newblock The theory of functions and sets of natural numbers, With a foreword
  by G. E. Sacks.

\bibitem{Odifreddi:1999}
Piergiorgio Odifreddi.
\newblock {\em Classical recursion theory. {V}ol. {II}}, volume 143 of {\em
  Studies in Logic and the Foundations of Mathematics}.
\newblock North-Holland Publishing Co., Amsterdam, 1999.

\bibitem{Sacks:1961*b}
Gerald~E. Sacks.
\newblock On suborderings of degrees of recursive unsolvability.
\newblock {\em Z. Math. Logik Grundlag. Math.}, 17:46--56, 1961.

\bibitem{Sacks:1963*b}
Gerald~E. Sacks.
\newblock On the degrees less than $\mbox{\bfseries 0}'$.
\newblock {\em Ann. of Math.}, 77:211--231, 1963.

\bibitem{Schnorr:1971*b}
C.-P. Schnorr.
\newblock A unified approach to the definition of random sequences.
\newblock {\em Math. Systems Theory}, 5:246--258, 1971.

\bibitem{Scott:1962}
Dana Scott.
\newblock Algebras of sets binumerable in complete extensions of arithmetic.
\newblock In {\em Recursive Function Theory}, volume~5 of {\em Proceedings of
  Symposia in Pure Mathematics}, pages 117--121, Providence, R.I., 1962.
  American Mathematical Society.

\bibitem{Shore:1978}
Richard~A. Shore.
\newblock On the $\forall\exists$-sentences of $\alpha$-recursion theory.
\newblock In R.~O.~Gandy J.~Fenstad and G.~E. Sacks, editors, {\em Generalized
  Recursion Theory II}, volume~94 of {\em Stud. Logic Foundations Math.}, pages
  331--354, Amsterdam, 1978. North--Holland Publishing Co.

\bibitem{Soare:1987}
Robert~I. Soare.
\newblock {\em Recursively Enumerable Sets and Degrees}.
\newblock Perspectives in Mathematical Logic, Omega Series. Springer--Verlag,
  Heidelberg, 1987.

\end{thebibliography}

\end{document}